\documentclass[12pt]{article}
\usepackage{amssymb}
\usepackage{graphicx}
\usepackage{amsfonts}
\usepackage{latexsym}
\usepackage{amsthm}
\usepackage{amsmath}

\usepackage{amssymb, amscd}

\usepackage{amsmath,amsfonts,graphicx,subfigure,url,amscd}

\usepackage{url}
\usepackage[T1]{fontenc}
\usepackage{color}
\usepackage{tabu}

\usepackage{bbm}

\numberwithin{equation}{section}

\newtheorem{theo}[equation]{Theorem}
\newtheorem{rem}[equation]{Remark}

\newtheorem{defin}[equation]{Definition}
\newtheorem{prop}[equation]{Proposition}
\newtheorem{cor}[equation]{Corollary}

\newtheorem{exam}[equation]{Example}

\begin{document}
\date{February 20, 2019}
 \title{{Combinatorics of unavoidable complexes\footnote{This research was supported by the Grants 174020 and
174034 of the Ministry of Education,
Science and Technological Development of the Republic of Serbia.}}}

\author{{Marija Jeli\'{c} Milutinovi\'{c}} \\
{\small Faculty of Mathematics}\\[-2mm]
{\small University of Belgrade}
\and{Du\v{s}ko Joji\'{c}} \\
 {\small Faculty of Science}\\[-2mm]
 {\small University of Banja Luka}
\and Marinko Timotijevi\'{c}\\
{\small Faculty of Science}\\[-2mm]
{\small University of Kragujevac}
\and Sini\v{s}a T. Vre\'{c}ica\\
{\small Faculty of Mathematics}\\[-2mm] {\small University of Belgrade}
\and Rade T. \v{Z}ivaljevi\'{c}\\
{\small Mathematical Institute}\\[-2mm] {\small SASA,
  Belgrade}\\[-2mm]}

\maketitle
\begin{abstract}\noindent
The partition number $\pi(K)$ of a simplicial complex $K\subseteq
2^{[n]}$ is the minimum integer $k$ such that for each partition
$A_1\uplus\ldots\uplus A_k = [n]$ of $[n]$ at least one of the
sets $A_i$ is in $K$. A complex $K$ is $r$-unavoidable if
$\pi(K)\leq r$. Motivated by the Van Kampen-Flores and
Tverberg type results, and inspired
by the `constraint method' \cite{bfz}, we study the
combinatorics  of $r$-unavoidable complexes. Emphasizing the interplay of ideas from combinatorial topology,  linear programming and fractional graph theory, we explore and compare extremal properties of examples arising in
topology (minimal triangulations) and combinatorics (hypergraph theory and Ramsey theory).
\end{abstract}

\section{Introduction}
\label{sec:intro}

The partition number $\pi(K)$ of a simplicial complex $K\subseteq 2^{[n]}$  is the minimum integer $k$ such that for each partition $A_1\uplus\ldots\uplus A_k = [n]$ of $[n]$ at least one of the
sets $A_i$ is in $K$. A simplicial complex $K$ is called {\em $r$-unavoidable} if $\pi(K)\leq r$.

A small partition number is important for applications, as illustrated by the inequality (Corollary~3.15 in \cite{jmvz}),
\begin{equation}\label{eqn:izvor-3}
{\rm Ind}_G (K^{\ast r}_\Delta) \geq n - \pi(K)
\end{equation}
where $r = \pi(K)=p^k$ is a prime power, $G=(\mathbb{Z}_p)^k$,
$K\subseteq 2^{[n]}$ is a simplicial complex, $K^{\ast r}_\Delta$ its $r$-fold deleted join, and ${\rm Ind}_G$ is an equivariant index function, see \cite{jmvz}. A variety of problems in topological combinatorics can be reduced to an appropriate index inequality (similar to \ref{eqn:izvor-3}), which explains why good upper bounds for $\pi(K)$ are relevant and interesting.

\medskip
The threshold characteristic $\rho(K)$  is the maximum real number $\alpha\geq 0$ such that for some probability measure  $\mu$ on $[n]$, the associated
 `threshold complex' $T_{\mu <\alpha}:=\{A\subset [n] \mid \mu(A)<\alpha\}$ is contained in $K$.
The fundamental relation between invariants $\pi(K)$ and $\rho(K)$ is the inequality
\begin{equation}\label{eqn:uvod-fund-ineq}
  \pi(K) \leq \lfloor 1/\rho(K)\rfloor +1
\end{equation}
(Proposition~\ref{prop:vazna-implikacija}) which, in essence, records the simple fact that for each partition $[n] = A_1\uplus\ldots\uplus A_k$, and each probability measure $\mu$ on $[n]$, there exists $i$ such that $\mu(A_i)\leq 1/k$ (the `pigeonhole principle' for measures).

\medskip

Our main objective is to construct {\em intrinsically non-linear} simplical complexes where the inequality (\ref{eqn:uvod-fund-ineq}) is strict and where the {\em non-linearity gap},
 \begin{equation}\label{eqn:gap-defin}
   \epsilon(K) :=    \lfloor 1/\rho(K)\rfloor +1 - \pi(K),
 \end{equation}
is as large as possible, relative to the size of $K$. In other words we search for examples of unavoidable complexes which are unavoidable for deeper reasons and cannot be detected by a simple application of the  pigeonhole principle (cf. \cite[Lemma~4.2]{bfz}).

\medskip

\begin{figure}[hbt]
\centering 
\includegraphics[scale=0.50]{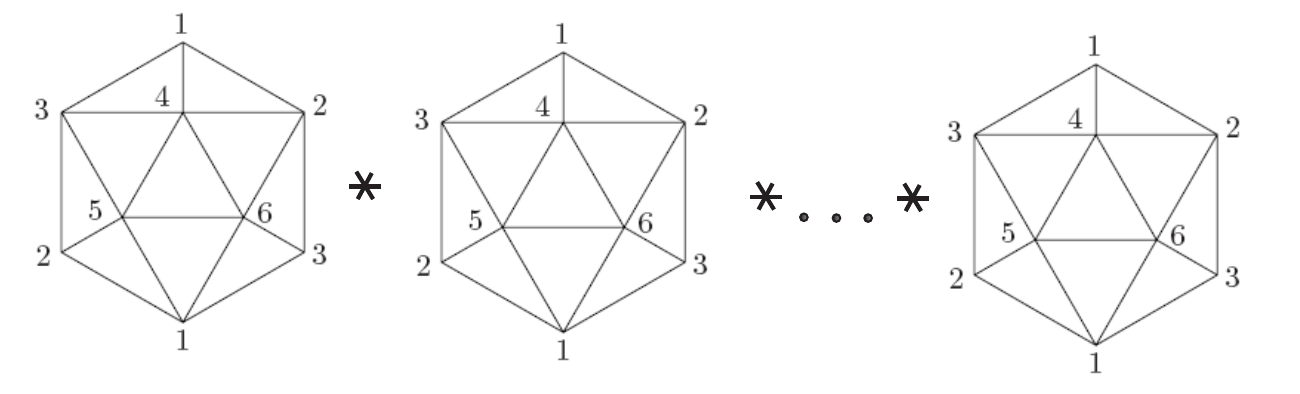}
 \caption{Join of $n$-copies of $\mathbb{R}P^2_6$.}
 \label{fig:poly-iko}
\end{figure}

Our examples are constructed as joins of self-dual (super-dual) complexes with a large group of automorphisms, as exemplified by the join $(\mathbb{R}P^2_6)^{\ast n}$ of $n$ copies of the $6$-vertex triangulation the real projective plane (Figure~\ref{fig:poly-iko}). Three main classes of self-dual (super-dual) complexes with a large automorphism group, analyzed in the paper, are the following:

 \begin{enumerate}
 \item[(1)] The triangulations of `projective planes' $\{[3], \mathbb{R}P^2_6, \mathbb{C}P^2_9, \mathbb{H}P^2_{15}\}$ with the minimum number of vertices, see \cite{bd-1, bk92, KB83, Lutz};
  \item[(2)] Complexes $\widehat{\mathcal{P}}_q$  associated to hypergraphs $\mathcal{P}_q$ of lines in a $q$-uniform finite projective plane, such as the Fano complex $\widehat{\mathcal{P}}_3$  (Section~\ref{sec:large-gap});

  \item[(3)]  Super-dual `Ramsey complexes'  $\mathfrak{R}_n$. For example $\mathfrak{R}_3 \subset 2^{E(K_6)} \cong 2^{[15]}$ is defined as the complex of all subgraphs $\Gamma\subset K_6$ of the complete graph $K_6$ on a set of $6$ vertices, such that the complement ${\Gamma}^c$ contains a triangle.
\end{enumerate}

The following table is a partial summary of our knowledge about the partition number $\pi$, the threshold characteristic $\rho$, and the non-linearity gap $\epsilon$ for joins of these complexes, as well as the quotient $\epsilon/w$ where $w$ is the  number of vertices.

\smallskip
\begin{equation}\label{eqn:tabela-nova}
  \begin{tabu}{c|c|c|c|c|c|c|}
  K & \widehat{\mathcal{P}}_3  & [3] & \mathbb{R}P^2_6 & \mathbb{C}P^2_9 & \mathbb{H}P^2_{15} & \mathfrak{R}_3\\\hline
  \pi(K^{\ast n}) & n+1 & n+1 & n+1 & n+1 & n+1 &  n+1\\ \hline
  \rho(K^{\ast n}) & 3/7n & 2/3n & 1/2n & 4/9n & 6/15n & \leq 6/15n\\ \hline
  \epsilon(K^{\ast n}) & \lfloor 4n/3\rfloor & \lfloor n/2\rfloor & n & \lfloor 5n/4\rfloor & \lfloor 3n/2\rfloor & ? \\ \hline
  w(K^{\ast n}) & 7n & 3n & 6n & 9n & 15n & 15n\\ \hline
  {\epsilon}/{w} & \approx 4/21 & \approx 1/6 & 1/6 & \approx 5/36 & \approx 1/10 &  ? \\\hline
  \end{tabu}
\end{equation}

\medskip\noindent
 On the basis of this information (and some other numerical and theoretical evidence), we conjecture that the inequality
 \begin{equation}\label{eqn:conjecture}
 \epsilon(K)/w(K) \leq 4/21 = \epsilon(\widehat{\mathcal{P}}_3^{\ast 3})/w(\widehat{\mathcal{P}}_3^{\ast 3})
 \end{equation}
 holds for all finite simplicial complexes $K$. The following theorem shows that the conjecture is true for flag simplicial complexes.

 \begin{theo}\label{thm:flag}
 Suppose that $K\subseteq 2^{[n]}$ is a `flag simplicial complex', i.e. the clique complex of a graph. Then,
  \begin{equation}\label{eqn:salient}
    \epsilon(K)/w(K) \leq 1/6 = \epsilon([3]^{\ast 2k})/w([3]^{\ast 2k}) \, .
 \end{equation}
 \end{theo}
A salient feature of the proof of Theorem~\ref{thm:flag}, which reduces the inequality (\ref{eqn:salient}) to a recent result of Choi, Kim, and O \cite{cko},  is the observation (conjectured by the anonymous referee!) that $1/\rho(K) +1$  (the right hand side of (\ref{eqn:uvod-fund-ineq})) can be interpreted as a fractional (linear programming) relaxation of the partition invariant $\pi(K)$.

\medskip

More explicitly, there is a link between the partition number $\pi(K)$ (respectively the threshold characteristic $\rho(K)$) and the matching number $\nu(F)$ (respectively the fractional matching
number $\nu_f(F)$) of the associated simple hypergraph $F\subseteq 2^{[n]}$ of minimal non-faces of $K$. In turn the inequality
(\ref{eqn:uvod-fund-ineq}) is reduced to $\nu(F) \leq \lfloor\nu_f(F)\rfloor$ and the study of the non-linearity gap (\ref{eqn:gap-defin}) is reduced to evaluating of the (integer part of) the fractional matching gap $\nu_f(F)- \nu(F)$.

\medskip
This connection allows us to obtain other results supporting the conjecture (\ref{eqn:conjecture}) as illustrated by the following theorem.

\begin{theo}\label{thm:uniform}
Let $q\geq 3$ and assume that $F\subset 2^{[n]}$ is a regular $q$-uniform hypergraph. Let $K\subseteq 2^{[n]}$ is the simplicial
complex associated to $F$, in the sense that $F$ is the collection of all minimal non-faces of $K$. Then
\begin{equation}\label{eqn:uniform}
\epsilon(K)/w(K) \leq \frac{(q-1)^2}{q(q^2-q+1)} \leq  4/21 \, .
\end{equation}

\end{theo}

\section{Unavoidable complexes and the constraint method}
\label{sec:part-inv-gl}

`Unavoidable complexes', originally introduced as `Tverberg unavoidable subcomplexes',  by Blagojevi\'{c}, Frick, and Ziegler \cite[Definition 4.1]{bfz}\footnote[2]{Some details in the original definition, such as the presence of a continuous map $f: \Delta_N \rightarrow \mathbb{R}^d$, are clearly superfluous. Indeed, the authors of \cite{bfz} added later that their main interest were the complexes large enough to be unavoidable for any continuous map.}, play the fundamental role in their `constraint method' (Gromov-Blagojevi\'{c}-Frick-Ziegler reduction).

\medskip
As an illustration how unavoidable complexes typically arise in applications, here we outline the basic idea of the `constraint method', as summarized in  \cite{Z17}.

\medskip
Suppose our goal is to prove a Tverberg-Van Kampen-Flores type result for a simplicial complex $K\subseteq \Delta^N$. More explicitly we want to show that for
each continuous map $f : K\rightarrow \mathbb{R}^d$ there exists vertex
disjoint simplices $\sigma_1,\ldots, \sigma_r\in K$ such that
$f(\sigma_1)\cap\ldots\cap f(\sigma_r)\neq\emptyset$.

\begin{equation}p
\begin{CD}
K @>f>> \mathbb{R}^d\\
@VeVV @ViVV\\
\Delta^N @>F>> \mathbb{R}^{d+1}
\end{CD}
\end{equation}
Let $K\subseteq \Delta^N$ be $r$-unavoidable. Assume that the continuous Tverberg theorem holds for the triple
$(\Delta^N, r, \mathbb{R}^{d+1})$, meaning that for each continuous
map $F: \Delta^N\rightarrow \mathbb{R}^{d+1}$ there exists a collection of
$r$ vertex disjoint faces $\Delta_1,\ldots, \Delta_r$ of
$\Delta^N$ such that $f(\Delta_1)\cap\ldots\cap
f(\Delta_r)\neq\emptyset$. This is the case, for example, if $r = p^k$ is a prime power and $N =
(r-1)(d+2)$.

Let $\bar{f}$ be an extension ($\bar{f}\circ e = f$) of
the map $f$ to $\Delta^N$. Suppose that $\rho :
\Delta^N\rightarrow \mathbb{R}$ is the function $\rho(x) := {\rm dist}(x,
K)$, measuring the distance of the point $x\in \Delta^N$ from $K$.
Define $F = (\bar{f}, \rho)  : \Delta^N\rightarrow \mathbb{R}^{d+1}$ and
assume that $\Delta_1,\dots, \Delta_r$ is the associated family of
vertex disjoint faces of $\Delta^N$, such that
$F(\Delta_1)\cap\ldots\cap F(\Delta_r)\neq\emptyset$. More
explicitly suppose that $x_i\in\Delta_i$ such that $F(x_i)=F(x_j)$
for each $i,j = 1,\ldots, r$. Since $K$ is $r$-unavoidable,
$\Delta_i\in K$ for some $i$. As a consequence $\rho(x_i)=0$, and
in turn $\rho(x_j)=0$ for each $j=1,\ldots, r$. If $\Delta_i'$ is
the minimal face of $\Delta^N$ containing $x_i$ then $\Delta_i'\in
K$ for each $i=1,\ldots, r$ and $f(\Delta_1')\cap\ldots\cap
f(\Delta_r')\neq\emptyset$.

For a more complete exposition and examples of
applications of the `constraint method' the reader is referred to
\cite{bfz}, see also \cite[Section~2.9(c)]{Grom-2} and \cite{L02}. An alternative approach, which relies on index inequalities similar to (\ref{eqn:izvor-3}), is developed in \cite{jvz-3}.

\begin{rem} {\rm
All `key examples' of $r$-unavoidable complexes, constructed and used in the seminal paper \cite{bfz} (Lemma~4.2), are {\em intrinsically linear} in the sense that they contain an $r$-unavoidable threshold complex, in which case the inequality (\ref{eqn:uvod-fund-ineq}) reduces to an equality.  Unavoidable threshold complexes provide at present the only general method for constructing simplicial complexes with a small partition number.  }
\end{rem}

\section{Threshold complexes and threshold characteristics}
\label{sec:threshold}

Let $\mu : 2^{[m]}\rightarrow \mathbb{R}_+$ be a monotone function ($A\subseteq B \Rightarrow
\mu(A)\leqslant \mu(B)$) such that $\mu(\emptyset)= 0$.
For a given `threshold' $\alpha\in \mathbb{R}_+$, the associated {\em threshold complex} is
  \begin{equation}\label{eqn:threshold-ini}
    T_{\mu \leq \alpha}:= \{A\in 2^{[m]} \mid
    \mu(A)\leq \alpha\}.
  \end{equation}

 From here on we focus on {\em linear threshold complexes} where $\mu = \mu_x$ is the measure (weight distribution)
 $\mu(A) := \sum_{i\in A} x_i$, associated to a non-negative weight vector $x\in \mathbb{R}^m_+$.

\begin{prop}\label{prop:concave}
If $\mu = \mu_x$ is a probability measure on $[m]$ ($\mu([n]) = \sum_{i\in [m]} x_i = 1$) then the associated threshold
complex $T_{\mu \leq 1/r}$ is $r$-unavoidable.
\end{prop}

\medskip\noindent
{\bf Proof:} Suppose that $[m] = A_1\uplus\ldots\uplus A_r$ is a
partition. Then at least one of the sets $A_i$ is in $T_{\mu  \leq
1/r}$, otherwise $1 =  \mu(A_1)+\ldots
+ \mu(A_r) >  1$.

\begin{defin}\label{def:linear}
An $r$-unavoidable simplicial complex  $K\subset 2^{[m]}$ is
called {\em linear} or {\em linearly realizable} if $K =
T_{\mu \leq 1/r}$ for some probability measure $\mu$ on $[m]$.
\end{defin}

Linearly realizable complexes are the simplest $r$-unavoidable
complexes and, in agreement with \cite[Section~4]{bfz}, they are also called
`pigeonhole complexes'. Note that a complex $K\subseteq 2^{[m]}\cong \{0,1\}^m\subset \mathbb{R}^m$ is linear (threshold) if and only if it can be separated by a hyperplane from its complement $2^{[m]}\setminus K \subset \{0,1\}^m$.

\medskip

It is interesting to study how large can be a threshold
complex which is contained in a simplicial complex $K$. We
introduce the {\em threshold characteristic} of the complex $K$ as
the unique number $\rho(K)$ such that
 \begin{equation}\label{eqn:char-ineq}
   \alpha< \rho(K) \quad
   \Longleftrightarrow \quad (\exists \mu\in\Delta_{[n]})~ T_{\mu \leq\alpha}\subset K \, ,
 \end{equation}
 where $\Delta_{[n]} = \{\mu_x \mid x\in \mathbb{R}^n_+ \mbox{ {\rm and} } \sum x_i = 1 \}$ is the simplex of probability measures on $[n]$.

 \medskip
 Note that $T_{\mu \leq\alpha} \subset K$ implies
 that $T_{\mu \leq\alpha +\epsilon} \subset K$
 for some $\epsilon> 0$, hence the set
 $\{\alpha \in \mathbb{R}_+ \mid T_{\mu \leq\alpha}
 \subset K\}$ is an open interval.

\begin{defin}\label{def:threshold-char}
  Let $K\subseteq 2^{[n]}$ be a simplicial complex and
  let $\Delta_{[n]}$ be the
  simplex of probability measures on $[n]$.
  The threshold characteristic $\rho(K)$ of $K$ is defined by
  \begin{align}
  \rho(K) &= \sup\{\alpha\in [0,+\infty] \mid
  (\exists \mu\in \Delta_{[n]})~ T_{\mu \leq\alpha}  \subseteq K\}
   \label{eq:rho-1}\\
    &=  \max\{\alpha\in [0,+\infty]
     \mid (\exists \mu\in \Delta_{[n]})~ T_{\mu <\alpha}
     \subseteq K\}\, .\label{eq:rho-2}
\end{align}
\end{defin}
 By definition $\rho(K) =
+\infty$ if and only if $K = 2^{[n]}$.

\begin{prop}\label{prop:vazna-implikacija}
 If $K\subseteq 2^{[n]}$ is a simplicial complex on $[n]$ then,
 \begin{equation}\label{eqn:vazna-nejednakost}
   \pi(K) \leq \lfloor {1}/{\rho(K)} \rfloor  +1.
 \end{equation}
\end{prop}

\medskip\noindent
{\bf Proof:} By definition if $1/r < \rho(K)$ then $T_{\mu \leq 1/r}\subset K$ for some probability
measure $\mu\in \Delta_{[n]}$. As a consequence of
Proposition~\ref{prop:concave} (the pigeonhole principle for
measures) we obtain the implication $1/r < \rho(K) \Longrightarrow
\pi(K)\leq r$, which implies the inequality (\ref{eqn:vazna-nejednakost}) if we choose
$r = \lfloor  1/\rho(K)\rfloor  + 1 $.
  \hfill $\square$

\medskip

The invariant $\rho(K)$ has a geometric interpretation in terms of {\em blocking polyhedra} $B(P)$ which were
 (for polyhedral sets $P$) introduced by Fulkerson, see \cite[Section~5.8]{sch03}.

\begin{defin}
  The convex set
  \begin{equation}\label{eqn:antidual}
    B(K) =
    \bigcap_{C\notin K}~\{ x\in \mathbb{R}^n_+ \mid \langle x,
    \chi_C\rangle \geq 1   \} \, ,
  \end{equation}
  is referred to as the {\em blocking polyhedron} of the simplicial complex $K\subseteq 2^{[n]}$.
\end{defin}

\begin{prop}\label{prop:thr-char-min}
Let  $\phi : \mathbb{R}^n\rightarrow \mathbb{R}$ be the
functional $\phi(x) = \langle x,
\mathbbm{1}\rangle =x_1+\dots + x_n$ and let
$\mathfrak{m}$ be the minimum value of $\phi$ on the
blocking polyhedron $B(K)$ of $K$.
Then, $$\rho(K) = 1/\mathfrak{m}\, .$$
\end{prop}

\medskip\noindent
{\bf Proof:} By Definition~\ref{def:threshold-char},
the threshold characteristic  $\rho(K)$ is the largest
$\alpha$ such that
\[
\Delta_{[n]}\cap \alpha B(K) \neq \emptyset \, ,
\]
or equivalently it is equal to the reciprocal of the
smallest $\mathfrak{m}$ such that
$
\mathfrak{m}\Delta_{[n]}\cap  B(K) \neq \emptyset
$.
It immediately follows  that $\mathfrak{m} =
\min\{\phi(x) \mid x\in B(K)\}$.
 \hfill $\square$

\medskip
The following proposition, useful in calculations,
is easily deduced from the observation that the
functional $\phi$ is $S_n$-invariant.

\begin{prop}\label{{{prop:vazna-simetrija-staro}}}
Let $G$ be the group of all permutations of $[n]$
which keep $K$ invariant and let $V = (\mathbb{R}^n)^G$ be
the associated invariant subspace of $\mathbb{R}^n$. Then,
\begin{equation}\label{eqn:simetrizacija}
\mathfrak{m} := {\min}\{\phi(x) \mid x\in B(K)\} =
{\min}\{\phi(x) \mid x\in B(K)\cap V\}.
\end{equation}
\end{prop}

The following direct consequence of (\ref{eq:rho-2})
is also very useful for explicit calculations of the
threshold characteristic $\rho(K)$.
(In agreement with Definition~\ref{def:threshold-char},
the minima and maxima are evaluated in the interval $[0,+\infty]$.)

\begin{prop}\label{prop:rho-calc}
  If $K\subset 2^{[n]}$ is a simplicial complex then,
 \begin{equation}\label{eq:rho-formula}
   \rho(K) =  \max_{\mu\in\Delta_{[n]}} \min_{C\notin K}\, \mu(C) \, .
 \end{equation}
    \end{prop}

\medskip
 Proposition~\ref{prop:rho-calc} can be
 considerably improved if $K$ admits a large group of symmetries.

\begin{prop}\label{prop:vazna-simetrija}
  Let $G$ be a group of all permutations of
  $[n]$ that keep the complex $K$ invariant.
  Let $\Delta_{[n]}^G \subset \Delta_{[n]}$ be the
  closed, convex set of all $G$-invariant probability measures.
  Then,
  \begin{equation}\label{eq:G-rho-formula}
   \rho(K) =  \max_{\mu\in\Delta_{[n]}^G}
   \min_{C\notin K}\, \mu(C) \, .
 \end{equation}
\end{prop}

\medskip\noindent
{\bf Proof:}  The proof follows from the observation
that for each $\mu\in \Delta_{[n]}$, the associated
$G$-average $\nu \in \Delta_{[n]}^G$, defined by
$\nu(C) := (1/\vert G\vert)\sum_{g\in G}~\mu(g(C))$,
satisfies the inequality
\[
\min_{C\notin K}~\mu(C)  \leq \min_{C\notin K}~\nu(C) \, .
\]

\medskip
The following important corollary records for
the future reference the simplest instance of
Proposition~\ref{prop:vazna-simetrija}.

\begin{cor}\label{cor:one-dim}
  Let $G$ be the group of all permutations of
  $[n]$ which keep the complex $K\subset 2^{[n]}$
  invariant. If the action of this group is transitive then,
  \begin{equation}\label{eqn:one-dim}
    \rho(K) = {\min}\left\{{\vert C\vert}/{n} \mid
    C\notin K \right\}.
  \end{equation}
\end{cor}

\medskip
In the following examples we collect some calculations
of the threshold characteristic which illustrate the
use of results from this section.

\begin{exam}\label{exam:full-symmetry} {\rm
If $K = {[n]\choose \leq k}\subset 2^{[n]}$ then $\rho(K) =
\frac{k+1}{n}$. Indeed, the automorphism group of $K$ is the full
symmetric group $S_n$ and we are allowed to apply
Corollary~\ref{cor:one-dim}. In the special case $k=0$ we obtain
that the threshold characteristic of the empty complex $\emptyset
\subset 2^{[n]}$ is $1/n$. If $K = {[n]\choose \leq k}$ and $n =
r(k+1)-1$ (for $r,k\geq 2$), then $\pi(K)=r$ and $\rho(K) =
\frac{k+1}{r(k+1)-1}$. This example shows that the  inequality
(\ref{eqn:vazna-nejednakost}) cannot be improved in general.
 Observe that $\rho(K)> 1 \Rightarrow \rho(K) = +\infty$
 so in particular $\rho(K) >  1$ if and only if $K = 2^{[n]}$.
 If $K = 2^{[n-1]}\subset 2^{[n]}$
 the group of automorphisms is $S_{n-1}$ and $\rho(K) = 1$,
 by a simple application of
 Proposition~\ref{prop:vazna-simetrija}.}
\end{exam}
\begin{exam}\label{example:cross}
{\rm Let $K = [2]^{\ast n} = [2] \ast \dots \ast [2] =
\partial (\diamondsuit_n)$ be the boundary complex of the
$n$-dimensional cross-polytope $\diamondsuit_n = {\rm Conv}\{e_i,
-e_i\}_{i=1}^n\subset \mathbb{R}^n$. In this case the group of
symmetries of $K$ is transitive on its vertices so $\rho(K)$ is
attained if $\mu$ is the uniform measure where $\mu(\pm e_i)=
1/2n$ for each $i\in [n]$. The minimal non-simplices  $\{e_i,
-e_i\}$ are all of the same cardinality which implies that
$\rho(K) = \mu(\{e_i, -e_i\}) = 1/n$.

This example shows that a complex $K$ can have a small
$\rho$-characteristic, and at the same time non-trivial homology
in the top dimension.

 }
\end{exam}

\section{Fractional relaxation of the partition invariant $\pi(K)$}
\label{sec:frac-relaxation}

\begin{theo}\label{thm:frac-relax}
Let $K\subseteq 2^S$ be a simplicial complex and $F \subseteq 2^S\setminus K$ the collection (hypergraph) of its minimal non-faces.
Then
\begin{equation}\label{thm:pocetna}
\pi(K) = \nu(F) + 1  \qquad  \rho(K) = \frac{1}{\nu_f(F)}
\end{equation}
where $\nu(F)$, respectively $\nu_f(F)$, are the matching and the fractional matching number of the hypergraph $F$.
\end{theo}

\begin{cor}  The inequality $\nu(F) \leq \nu_f(F)$ is equivalent to the inequality (\ref{eqn:vazna-nejednakost}).
In other words the fundamental inequality (\ref{eqn:uvod-fund-ineq}) is a consequence of the LP relaxation inequality for the
matching number $\nu(F)$.
\end{cor}

Recall that the matching number $\nu(F)$ of a hypergraph $F\subseteq 2^S$ is the cardinality $\nu$ of a largest subfamily
$F'\subseteq F$ such that $A\cap B = \emptyset$ for each pair $A, B$ of distinct elements in $F'$. This number can be also
described as the solution  of a maximization problem over non-negative integers (integer programming).

\medskip
Assuming $\vert S\vert = n $ and $\vert F\vert = m$, let $M$ be the $n\times m$ {\em incidence matrix} of the hypergraph $F\subseteq 2^S$. We denote by $\mathbbm{1} = (1,1,\dots, 1)^T$
the corresponding column vector in both $\mathbb{R}^n$ and $\mathbb{R}^m$. Then, by interpreting $x = \chi_{F'}\in \mathbb{Z}^m$ as the characteristic
function of the family $F'\subseteq F$, we observe that the matching number $\nu(F)$ is the solution of the following integer programming problem:
\begin{equation}\label{eqn:integer}
\begin{array}{l}
\max \,\,  \mathbbm{1}^T  \cdot x \, , x\in \mathbb{Z}^m   \\
Mx  \leq \mathbbm{1}   \\
x  \geq 0
\end{array}
\end{equation}
The LP-relaxation of (\ref{eqn:integer}) is the linear program obtained by replacing the condition $x\in \mathbb{Z}^m$ by the condition $x\in \mathbb{R}_+^m$ (equivalently $x\in [0,1]^m$).
This linear program is exhibited in (\ref{eqn:LP-duality}) together with the corresponding dual program.
\begin{equation}\label{eqn:LP-duality}
\begin{array}{cccccccc}
\max \,\,  \mathbbm{1}^T  \cdot x  &&&&&&&  \min \,\,  y^T\cdot \mathbbm{1}  \\
Mx  \leq \mathbbm{1}  &&&&&&&  y^T M  \geq \mathbbm{1} \\
x  \geq 0   &&&&&&& y  \geq 0
\end{array}
\end{equation}
By definition $\nu_f(F)$, the fractional matching number, is the solution to the maximization problem (\ref{eqn:LP-duality}) (on the left) which is, as a consequence of LP-duality, the same as the solution to the minimization problem (\ref{eqn:LP-duality}) (on the right). The inequality $\nu(F)\leq \nu_f(F)$ is an immediate consequence and the vanishing (and non-vanishing!) of the associated {\em ``fractional duality gap''}  $\nu_f(F) - \nu(F)$ is an interesting and well studied question.

\bigskip\noindent
{\bf Proof of Theorem~\ref{thm:frac-relax}:}  The equality $\pi(K) = \nu(F) + 1$ is straightforward. Indeed, $k< \pi(K)$ is equivalent to the existence of a partition
$A_1\sqcup\dots\sqcup A_k = S$ where $A_i\notin K$ for each $i=1,\dots, k$.

\medskip
For the second equality we interpret the solution $y_0$ of the dual linear program in (\ref{eqn:LP-duality}) as the weight vector of a measure
$\mu_{y_0}$ on $S\cong [n]$. Then $\mathfrak{m}:= \mu_{y_0}([n]) = \langle y_0, \mathbbm{1} \rangle$ and $\mu_{y_0}(C)\geq 1$ for each $C\in F = K^c$.
Moreover $\mathfrak{m}$ is the smallest value that can arise from a feasible vector in the minimization
problem (\ref{eqn:LP-duality}). It follows that Proposition~\ref{prop:thr-char-min} is precisely the statement needed to complete the proof.   \hfill $\square$

\section{Intrinsically non-linear unavoidable complexes}
\label{sec:intrinsic}

\begin{defin}\label{def:intrin} A simplicial complex $K\subseteq 2^S$ is
intrinsically linear if the inequality
(\ref{eqn:vazna-nejednakost}) reduces to an equality. The {\em
non-linearity gap} of a complex $K$ is
 \begin{equation}\label{eqn:gap-defin-bis}
   \epsilon(K) :=    \lfloor 1/\rho(K)\rfloor +1 - \pi(K),
 \end{equation}
 and $K$ is intrinsically non-linear if $\epsilon(K)>0$. The relative non-linearity gap of $K$ is $\bar{\epsilon}(K) = \epsilon(K)/w(K)$ where $w(K)$ is the cardinality of $S$.
\end{defin}

\noindent
The following remark clarifies why the complexes with
$\epsilon(K)=0$ are called intrinsically linear.

\begin{rem}\label{prop:kad-je-jednakost}
{\rm  Let $K\subseteq 2^{[n]}$ be a simplicial complex such
  that $\pi(K) = r$.
  Then $\pi(K) = \lfloor{1}/{\rho(K)}\rfloor + 1$
  if and only if  $T_{\mu \leq 1/r}\subseteq K$ for
  some probability measure $\mu\in\Delta_{[n]}$.
  In other words if $\pi(K) = r$ then $K$ is an
  intrinsically linear complex if and only if it contains a
  linear $r$-unavoidable complex $T_{\mu \leq 1/r}$.}
\end{rem}

By the hypergraph description of $\pi(K)$ and $\rho(K)$ (Section~\ref{sec:frac-relaxation}) there is an equality
$\epsilon(K) = \lfloor\nu_f(F)\rfloor - \nu(F)$.  This observation focuses our attention to the fractional matching gap $\epsilon_f(K) = \nu_f(F) - \nu(F)$ and its relative version $\bar{\epsilon}_f(K) = (\nu_f(F) - \nu(F))/w(K)$.

\begin{prop}\label{prop:convex}
Let  $K_1\subseteq 2^{S_1}$ and $K_2\subseteq 2^{S_2}$ be a pair of simplicial complexes and let $F_1\subseteq 2^{S_1}$ and  $F_2\subseteq 2^{S_2}$ be the associated hypergraphs of minimal non-faces (Section~\ref{sec:frac-relaxation}).
Then the hypergraph associated to the join $K = K_1\ast K_2 \subseteq 2^{S}$ (where $S = S_1\uplus S_2$) is the disjoint union
$F = F_1\uplus F_2 \subseteq 2^S$. Moreover,
\[
\bar{\epsilon}_f(K_1\ast K_2) = \frac{n_1}{n_1+n_2}\bar{\epsilon}_f(K_1) + \frac{n_2}{n_1+n_2}\bar{\epsilon}_f(K_2)
\]
where $n_1 = w(K_1) = \vert S_1\vert$ and $n_2 = w(K_2) = \vert S_2\vert$.
\end{prop}

\medskip\noindent
{\bf Proof:} Both $\nu$ and $\nu_f$ are additive with respect to the disjoint sum of hypergraphs. For example a fractional matching $\phi : F_1\uplus F_2 \rightarrow \mathbb{R}_+$ is simply a sum of two fractional matchings on $F_1$ and $F_2$, etc. \hfill $\square$

\begin{cor}\label{cor:rel-gap-joins}
Let $K^{\ast n}$ be the join of $n$ copies of $K\subseteq 2^S$ and let $F\subseteq 2^S$ be the hypergraph associated to $K$. Then
\begin{equation}\label{eqn:gap-joins}
{\epsilon}_f(K^{\ast n}) = n \epsilon_f(K) \qquad  \mbox{ {\rm and} }  \qquad
\bar{\epsilon}_f(K) = \bar{\epsilon}_f(K^{\ast n}) =  \bar{\epsilon}(K^{\ast n})
\end{equation}
where the rightmost equality holds under condition that $n\nu_f(F)$ is an integer.
\end{cor}

\section{Complexes with large non-linearity gap}
 \label{sec:large-gap}

By Corollary~\ref{cor:rel-gap-joins} the non-linearity gap ${\epsilon}(K^{\ast n})$  behaves essentially as a linear function in $n$, with the rate of increase $\bar{\epsilon}_f(K)$. Moreover, as a consequence of Proposition~\ref{prop:convex},
\begin{equation}\label{eqn:conv-ineq}
\bar{\epsilon}_f(K_1\ast K_2) \leq \max \{\bar{\epsilon}_f(K_1), \bar{\epsilon}_f(K_2)\} \, .
\end{equation}
 For this reason we focus our attention to examples of $2$-unavoidable complexes ($\pi(K) = 2$) where $\bar{\epsilon}_f(K)$ is as large as possible. Note that a simplicial complex $K$ is $2$-unavoidable if and only if the associated hypergraph is intersecting ($A\cap B\neq\emptyset$ for each $A, B\in F$).

\begin{enumerate}
\item
Let  $\mathbb{F}$ be a finite field of order $p^k$ and $P(\mathbb{F}^3)$ the corresponding finite projective plane. Let
${\mathcal{P}}_q$ be the $q$-uniform hypergraph of all lines in $P(\mathbb{F}^3)$ ($q = p^k+1)$ and let
 $\widehat{\mathcal{P}}_q$ the associated simplicial complex ($A\in \widehat{\mathcal{P}}_q$ if and only if $A$ does not contain a line). For more information on ${\mathcal{P}}_q$ the reader is refereed to \cite{b89, f81}.

\item
From the list of $11$ exceptional, vertex-minimal
triangulations of manifolds exhibited in \cite[Table~2]{Lutz} we select the
`projective planes' $\{[3], \mathbb{R}P^2_6, \mathbb{C}P^2_9,
\mathbb{H}P^2_{15}\}$, since they are all Alexander self-dual complexes
($\pi(K) = 2$) and have a vertex-transitive group of symmetry
(which makes the $\rho$-invariant easily computable). More
information about these important complexes can be found in \cite{KB83, bd-1, bk92, Go}.

\item   Let $\mathfrak{r} = \mathfrak{r}(n)$ be the Ramsey number (the minimum number $m$ such that each
subgraph of $K_m$ contains either a clique of size $n$ or an independent
subset of size $n$). Let $V = {{[\mathfrak{r}]}\choose{2}}$ be  the set of edges of the complete
graph $K_\mathfrak{r}$ on $\mathfrak{r}$ vertices. Let
$\mathfrak{R}_n\subset 2^V$ be the collection of all subsets
(graphs)  $\Gamma \subset K_\mathfrak{r}$ such that the complement
${\Gamma}^c = V\setminus \Gamma$ contains a copy of the complete
graph $K_n$ with $n$ vertices. By Ramsey's theorem
$\mathfrak{R}_n$ is a $2$-unavoidable simplicial complex.
\end{enumerate}

Following the procedures described in Section~\ref{sec:threshold}
(especially Proposition~\ref{prop:vazna-simetrija} and
Corollary~\ref{cor:one-dim}), we calculate the threshold characteristic and other invariants of these complexes.

\medskip\noindent
$\widehat{\mathcal{P}}_q$

\smallskip\noindent
Obviously $\pi(\widehat{\mathcal{P}}_q)=2$. Either by the results of Section~\ref{sec:threshold} or by Section~\ref{sec:frac-relaxation} and \cite[Section 1.3]{f81}, we have $\rho(\widehat{\mathcal{P}}_q)=q/(q^2-q+1)$. It follows that
in this case $\epsilon_f(\widehat{\mathcal{P}}_q) =  (q-1)^2/q$  and $\bar{\epsilon}_f(\widehat{\mathcal{P}}_q) = (q-1)^2/q(q^2-q+1)$.
In the special case of the `Fano complex' $\widehat{\mathcal{P}}_3$ (corresponding to the seven element Fano plane ${\mathcal{P}}_3$), one obtains $\bar{\epsilon}_f(\widehat{\mathcal{P}}_3) = 4/21$.

\medskip\noindent
$\mathbb{R}P^2_6$

\smallskip\noindent
The minimum triangulation of the real projective plane has six
vertices  (Figure~\ref{fig:poly-iko}). The group of simplicial
automorphisms of $\mathbb{R}P^2_6$ is vertex-transitive. The size
of the smallest non-face is $3$. It follows that
$\rho(\mathbb{R}P^2_6)= 1/2$.

\medskip\noindent
$\mathbb{C}P^2_9$

\smallskip\noindent
The minimum triangulation of the complex projective plane has $9$
elements. The size of the smallest non-face is $4$. The group
of simplicial automorphisms of $\mathbb{C}P^2_9$ is also
vertex-transitive. This can be easily deduced from the description
of this group given in \cite{KB83} (see also
\cite[Section~2]{bd-1}). It follows that $\rho(\mathbb{C}P^2_6)=
4/9$.

\medskip\noindent
$\mathbb{H}P^2 = M^8_{15} = \mathbb{H}P^2_{15}$

\smallskip\noindent
Brehm and K\"{u}hnel constructed in \cite{bk92} $PL$-isomorphic
simplicial complexes $M^8_{15}, \widetilde{M}^8_{15},
{\widetilde{\widetilde{M}^8_{15}}}$ on $15$ vertices  and
conjectured that they triangulate the quaternionic projective
plane $\mathbb{H}P^2$ (this conjecture was recently confirmed by
Gorodkov in \cite{Go}). The most symmetric among them is
$M^8_{15}$, which is invariant under  vertex-transitive action of
the group $A_5$ (see the Theorem on page 169 in \cite{bk92}).
Moreover the cardinality of the minimum size non-face in
$M^8_{16}$ is $6$. From here we deduce that $\rho(M^8_{15}) =
6/15$.

\medskip\noindent
$\mathfrak{R}_3$

\begin{figure}[hbt]
\centering
\includegraphics[scale=0.35]{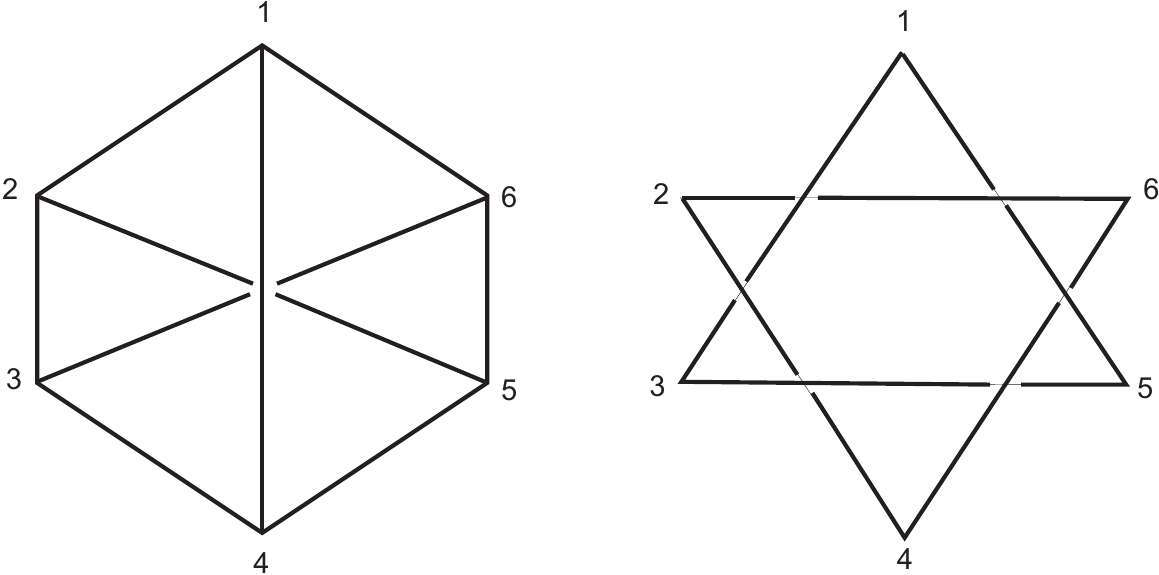}
\caption{A maximal graph without triangles and its
complement.}\label{fig:Ramsey-1}
\end{figure}

\smallskip\noindent
Let $\mathfrak{R}_3$ be the complex of all graphs $\Gamma \subset K_6$ such that the complement ${\Gamma}^c = V\setminus \Gamma$ contains a triangle $K_3$. Then $\pi(\mathfrak{R}_3) = 2, \rho(\mathfrak{R}_3) = 6/15$ and $w(\mathfrak{R}_3) = 15$  (Figure~\ref{fig:Ramsey-1}).

  More generally, define $K\subset 2^E$ as a simplicial complex on the set $E={[m]\choose 2}$ of all edges in the complete graph $K_m$, where $m\geq \mathfrak{r}(n+1)$ and $m = tn$ is divisible by $n$. By definition $\Gamma \in K$ if and only if $\Gamma$, interpreted as graph on $[m]$, contains an independent set of size $n+1$.  Then, $\pi(K) = 2,  \rho(K) \leq {1}/{n}$ and $\nu(K) =  {nt\choose 2}$.

\bigskip
The following theorem summarizes the results of calculations collected in Table~(\ref{eqn:tabela-nova}).

\begin{theo}\label{thm:vazna-2} Let $K$ be one of the
following complexes, $\{\widehat{\mathcal{P}}_q, [3], \mathbb{R}P^2_6,
\mathbb{C}P^2_9, \mathbb{H}P^2_{15}, \mathfrak{R}_3\}$.
 Let $K^{\ast n}$ be the join of $n$ copies of $K$.
 For a given complex $L\subseteq 2^S$, let $\pi(L), \rho(L),
 \epsilon(L), w(L)$ be respectively, the partition invariant,
 the threshold characteristic, the non-linearity gap,
 and the number of vertices of $L$.
 Then the  values of these invariants for the
 complex $L = K^{\ast n}$ are collected in
 Table~(\ref{eqn:tabela-nova}).
\end{theo}

\bigskip\noindent
{\bf Proof:}  All complexes $K$ listed in the theorem are $2$-unavoidable. The first five of them are actually
Alexander self-dual \cite{KB83, bd-1, bk92} (which means that they are minimal $2$-unavoidable). In both cases $\pi(K)=2$.

\smallskip
The second row in  Table~(\ref{eqn:tabela-nova}) (evaluating the invariant $\rho$) is
computed  by the methods developed in Section~\ref{sec:threshold}
which uses Proposition~\ref{prop:vazna-simetrija} as the key tool. For $K^{\ast n}$ one uses the additivity of invariants $\nu$ and $\nu_f$ and their link with $\pi$ and $\rho$ (Section~\ref{sec:frac-relaxation}).
The rest of Table~\ref{eqn:tabela-nova} is completed by simple calculation. \hfill $\square$

\subsection{Theorems~\ref{thm:flag} and \ref{thm:uniform}, and the origin of Conjecture~(\ref{eqn:conjecture})}

\noindent
{\bf Proof of Theorem~\ref{thm:flag}:} By a recent result of Choi, Kim, and O (see \cite[Corollary 7]{cko}),
for any $n$-vertex graph $G$ there is an inequality $\nu_f(G) - \nu(G)\leq n/6$, with the equality only when $G$ is a disjoint union of triangles $K_3$. This, together with the results from Section~\ref{sec:frac-relaxation}, immediately leads to the inequality (\ref{eqn:salient}). \hfill $\square$

\medskip\noindent
{\bf Proof of Theorem~\ref{thm:uniform}:} By \cite{b89} (Corollary 2 on p. 104) if $F\subseteq 2^{[n]}$ is a regular $q$-uniform hypergraph ($q\geq 3$), then $\nu \geq \frac{n}{q^2-q+1}$. Together with the equality $\nu_f(F) = \frac{n}{r}$, which holds for all $q$-uniform hypergraphs (\cite{b89}, Theorem 7 on p. 94)), this leads to the inequality
\begin{equation}\label{eqn:uniform-pf}
  \frac{\nu_f(F) - \nu(F)}{n} \leq \frac{(q-1)^2}{q(q^2-q+1)}
\end{equation}
which is in light of Theorem~\ref{thm:frac-relax} equivalent to (\ref{eqn:uniform}). \hfill $\square$

\medskip
The largest value of the right hand side of (\ref{eqn:uniform-pf}) is attained for $q=3$. This corresponds to the case of the $3$-uniform Fano  plane $\mathcal{P}_3$ when the inequality (\ref{eqn:uniform-pf}) reduces to the equality $\frac{1}{7}(\nu_f(\mathcal{P}_3) - \nu(\mathcal{P}_3)) = \frac{4}{21}$. It follows that in the conjectured inequality (\ref{eqn:conjecture}) the equality is attained for all complexes $K = \widehat{\mathcal{P}}_3^{\ast 3k}$, where $k\geq 1$.

 \bigskip\noindent
 {\bf Acknowledgements:} We acknowledge kind
 remarks by P. Landweber, R. Meshulam and S. Todor\v cevi\'{c}, at an earlier stage of the project, and of G. Simonyi who emphasized the importance of references \cite{cko,f81}. Two of the authors
 (M.~Jeli\'{c} and R.~\v{Z}ivaljevi\'{c})  participated in the program
 `Topology in Motion' at the Institute for Computational and
 Experimental Research in Mathematics (ICERM, Brown University).
 With great pleasure they acknowledge the support, hospitality and
 excellent working condition at ICERM. We are much obliged to the
 referees for numerous remarks and suggestions which considerably
 improved the presentation of results in the paper. In particular we are indebted to the referee who conjectured that $1/\rho(K) +1$ is a fractional relaxation of $\pi(K)$.


\begin{thebibliography}{00000}

\bibitem[BD94]{bd-1} B.~Bagchi, B.~Datta. On
$9$-vertex complex projective plane. Geom. Dedicata 50 (1994), 1--13.


\bibitem[B89]{b89} C. Berge. \textit{Hypergraphs: Combinatorics of Finite Sets}, North-Holland (1989).


\bibitem[BFZ]{bfz}
P.V.M.~Blagojevi{\'c}, F.~Frick, G.M.~Ziegler.
\newblock Tverberg plus constraints.
\newblock {\em Bull. Lond. Math. Soc.}, 46:953--967, 2014.


 \bibitem[BK92]{bk92} U.~Brehm, W.~K\"{u}hnel. 15-vertex
 triangulations of an 8-manifold.  {\em Math. Ann.},
 294, Issue 1,   167--193 (1992).

\bibitem[CKO]{cko} I. Choi, J. Kim, S.O.
The difference and ratio of the fractional matching number and the matching number of graphs,
\textit{Discrete Math.} 339 (2016), 1382--1386.


\bibitem[F81]{f81} Z.~F\" uredi. Maximum degrees and fractional matchings in uniform hypergraphs, \textit{Combinatorica} 1 (1981)
155--162.


\bibitem[Go]{Go} D.~Gorodkov. A 15-vertex
triangulation of the quaternionic projective plane.
arXiv:1603.05541 [math.AT]

\bibitem[Gr10]{Grom-2} M.~Gromov. Singularities, expanders
and topology of
maps. Part 2: From combinatorics to topology via algebraic
isoperimetry  {\em Geom. Funct. Anal.}  20 (2010), 416--526.





\bibitem[JVZ]{jvz-3}
D.~Joji{\'c}, S.T.~Vre{\'c}ica, R.T.~{\v Z}ivaljevi{\'c}.
\newblock Topology and combinatorics of
`unavoidable complexes', arXiv:1603.08472v1 [math.AT], (unpublished prepreint).

\bibitem[JMVZ]{jmvz}
D.~Joji{\'c}, W.~Marzantowicz, S.T.~Vre{\'c}ica, R.T.~{\v Z}ivaljevi{\'c}.
\newblock Topology of unavoidable complexes, arXiv:1603.08472 [math.AT].

\bibitem[KB83]{KB83} W.~K\"{u}hnel,
T.F.~Banchoff. The $9$-vertex complex projective plane.
{\em Math. Intell.} 5(3), 11--22 (1983).

\bibitem[Lon02]{L02}
M. de Longueville.
Notes on the topological Tverberg theorem,
\textit{Discrete Math.}, 241:207–-233, 2001.
Erratum: \textit{Discrete Math.}, 247:271–-297, 2002.

\bibitem[L]{Lutz} F.H. Lutz. Triangulated
manifolds with few vertices: Combinatorial
manifolds. 2005. Preprint, 2005; arXiv:math/0506372.




\bibitem[Sch03]{sch03} A. Schrijver.
\textit{Combinatorial Optimization:
Polyhedra and Efficiency}. Springer, Berlin 2003.

\bibitem[\v Z17]{Z17}
R.T.~\v Zivaljevi\'{c}. Topological methods in discrete geometry.
Chapter 21 in \textit{Handbook of Discrete and Computational
Geometry}, third ed., J.E. Goodman, J. O'Rourke, and C.D. T\'{o}th,
CRC Press LLC, Boca Raton, FL, 2017.






\end{thebibliography}
\end{document}